\documentclass[10pt,twoside]{article}
\usepackage{amssymb}
\usepackage{latexsym}
\newtheorem{theorem}{Theorem}[section]

\newtheorem{lemma}{Lemma}[section]
\newtheorem{proposition}{Proposition}[section]
\newtheorem{remark}{Remark}[section]

\def\proof{\mbox {\it Proof.~}}
\makeatletter\def\theequation{\arabic{section}.\arabic{equation}}\makeatother

\begin{document}
\title{
\vspace{0.5in}
{\bf\Large  Nonlinear Boundary Value Problems via Minimization on Orlicz-Sobolev Spaces}\hspace{2mm}}
\author{{\bf\large J. V.  Goncalves}\footnote{The author  acknowledges the support by CNPq/Brasil. }\hspace{2mm}\\
{\it\small  Instituto de Matem\'atica e Estat\'istica},
{\it\small Universidade Federal de Goi\'as}\\
{\it\small  74001-970 Goi\^ania, GO - Brasil}\\
{\it\small e-mail: goncalves.jva@pq.cnpq.br}\\
\vspace{1mm}\\
{\bf\large M. L. M. Carvalho}\footnote{The author acknowledges the support by CAPES/Brasil.}\hspace{2mm}\\
{\it\small Departamento de Matem\'atica},
{\it\small Universidade Federal de Goi\'as}\\
{\it\small  75804-020  Jata\'i, GO - Brasil}\\
{\it\small e-mail: marcosleandro@jatai.ufg.br}\vspace{1mm}}
\maketitle
\begin{center}
{\bf\small Abstract}
\vspace{3mm}
\hspace{.05in}\parbox{4.5in}
{{\small  We develop arguments on convexity and minimization of  energy functionals  on Orlicz-Sobolev spaces to investigate existence of solution to the equation
$\displaystyle  -\mbox{div} ( \phi(|\nabla u|) \nabla u) = f(x,u) + h~ \mbox{in}~ \Omega$
 under Dirichlet boundary conditions, where  $\Omega \subset {\bf R}^{N}$  is a bounded smooth domain,  $\phi : (0,\infty)\longrightarrow (0,\infty)$  is a suitable continuous function and   $f: \Omega \times {\bf R} \to {\bf R}$ satisfies the  Carath\'eodory conditions,  while $h: \Omega \to {\bf R}$ is a  measurable function.
}}
\end{center}
\noindent
{\it \footnotesize 2000 Mathematics Subject Classifications}. {\scriptsize 35J25, 35J60}.\\
{\it \footnotesize Key words}. {\scriptsize quasilinear equations, minimization, Orlicz-Sobolev spaces.}

\section{\bf Introduction}
\def\theequation{1.\arabic{equation}}\makeatother
\setcounter{equation}{0}

We develop arguments on convexity and minimization of  energy functionals  on Orlicz-Sobolev spaces to investigate existence of solution to the problem
\begin{equation}\label{prob}
\displaystyle  -\mbox{div} ( \phi(|\nabla u|) \nabla u) = f(x,u) + h~~ \mbox{in}~
\Omega,~~ u = 0~  \mbox{on}~ \partial \Omega,
\end{equation}
 where  $\Omega \subset {\bf R}^{N}$  is a bounded smooth domain and $\phi : (0,\infty)\longrightarrow (0,\infty)$  is a continuous function satisfying
\vspace{1mm}

\begin{itemize}
  \item[($\phi_1$)] $\mbox{(i)}~\displaystyle  \lim_{s\rightarrow 0}s\phi (s) =0$,~~ $\mbox{(ii)}~\displaystyle  \lim_{s\rightarrow \infty} s\phi (s) = \infty$,~
  \item[($\phi_2$)] $s \mapsto s\phi (s)$ is nondecreasing in $(0, \infty),$
\end{itemize}
\noindent We extend $s \mapsto s\phi (s)$ to $\mbox{\bf R}$ as an odd function and consider the associated even potential
$$
\Phi(t):=\int_0^t s\phi(s)ds,~t\in {\bf R}.
$$
\noindent It follows from the continuity of  $\phi$, $(\phi_1)$ and $(\phi_2)$ that $\Phi$ is  increasing  and convex.
The function $f: \Omega \times {\bf R} \to {\bf R}$ satisfies the  Carath\'eodory conditions,  while $h: \Omega \to {\bf R}$ is assumed  measurable. Further conditions will be imposed upon $f$ and $h$ in a while.

\noindent  Next we will introduce some notations concerning Orlicz and Orliccz-Sobolev spaces. The function  $\Phi$ is said to satisfy the  $\Delta_2$-condition, $\Phi \in \Delta_2$ for short, if
$$
   \mbox{there is a constant}  K>0~\mbox{such that }~ \Phi(2t)\leq K\Phi(t),~~t\geq 0.
$$
\noindent The complementary function $\widetilde{\Phi}$ associated to $\Phi$ is defined by
$$
\widetilde{\Phi}(t):=\displaystyle\max_{s\geq 0}\{st-\Phi(s)\},~t\geq 0.
$$
We recall, (cf.  \cite[thm 3, pg. 22]{Rao1}), that  $\Phi,\widetilde{\Phi} \in \Delta_2$ iff there are  $ \ell, m\in(1,\infty)$ such that
\begin{equation}\label{cond_phi}
 \ell\leq\frac{t^2\phi(t)}{\Phi(t)}\leq m, ~~t> 0,
\end{equation}
 We shall assume from now on that both $\Phi$ and $\widetilde{\Phi}$ satisfy the $\Delta_2$-condition.

\noindent We recall, see e.g. Adams \& Fournier \cite{A}, that the Orlicz Space associated with $\Phi$ is given by
$$
\displaystyle L_\Phi(\Omega):=\left\{u:\Omega\longrightarrow{\mbox{\bf R}}~ \mbox{measurable}~| ~ \int_\Omega \Phi\left(\frac{u(x)}{\lambda}\right)<+\infty~ \mbox{for some}~\lambda>0\right\}.
$$
\noindent It is known, (cf. \cite{A} ), that  the expression
$$\|u\|_\Phi=\inf\left\{\lambda>0~|~\int_\Omega \Phi\left(\frac{u(x)}{\lambda}\right)\leq 1\right\}$$
\noindent  defines a norm  in $L_{\Phi}(\Omega)$ named Luxemburg norm. By \cite[lemma D2]{Clement-1},
\begin{equation}\label{emb-ell}
L_\Phi(\Omega)\stackrel{\hookrightarrow}{\mbox{\tiny cont}} L^{\ell}(\Omega).
\end{equation}

\noindent  The corresponding Orlicz-Sobolev space, (also denoted   $W^1L_\Phi(\Omega)$),  is defined as
$$
 W^{1, \Phi}(\Omega) =\Big\{u \in L_\Phi(\Omega)~|~ \frac{\partial u}{\partial x_i} \in L_\Phi(\Omega),~ i=1,...,N \Big\}
$$
\noindent The usual Orlicz-Sobolev norm of $ W^{1, \Phi}(\Omega)$ is
$$
  \displaystyle \|u\|_{1,\Phi}=\|u\|_\Phi+\sum_{i=1}^N\left\|\frac{\partial u}{\partial x_i}\right\|_\Phi.
$$
\noindent Since we are assuming that $\Phi$ and $\widetilde{\Phi}$ satisfy the $\Delta_2$-condition,  $L_{\Phi}(\Omega)$  and $W^{1,\Phi}(\Omega)$  are separable, reflexive,  Banach spaces, see e.g. \cite{A}. We also set
$$
W_0^{1,\Phi}(\Omega)=\overline{{\bf C}_0^\infty(\Omega)}^{W^{1,\Phi}(\Omega)}.
$$
\noindent One shows that $u\in W_0^{1,\Phi}(\Omega)$ means that~ $\displaystyle u = 0~~\mbox{on}~~ {\partial\Omega}$~ in the trace sense,  cf. Gossez \cite{Gz1}.
\vspace{1mm}

\noindent In order to state our main results consider the potential function of $f$,
$$
\displaystyle F(x,t):=\int_0^t f(x,s)ds
$$
\noindent and the limit
$$
\displaystyle A_\infty(x):=\displaystyle \limsup_{|s|\rightarrow\infty}\frac{F(x,s)}{|s|^\ell}.
$$
\noindent We shall assume that that there exist a number $ A \geq 0$ and a nonnegative function  $B \in L^1(\Omega)$ such that
\begin{equation}\label{cond_F}
        \displaystyle
       F(x,s) \leq A|s|^{\ell} + B(x),~~  s \in \mbox{\bf R},~~ \mbox{a.e.}~~ x \in {\Omega}.
         \end{equation}
\noindent  From now on suppose $1<\ell,m<N$.

\noindent We set
$$
\displaystyle\Phi_*^{-1}(t) := \int_0^t\frac{\Phi^{-1}(s)}{s^\frac{N+1}{N}}ds,~~t > 0.
$$
\noindent The {\rm critical exponent}  function of $\Phi$,  $\Phi_*$, is defined  as the inverse function of $\Phi_*^{-1}$.
It is known  that  $\Phi_*$ is an N-function, see e.g. Donaldson \& Trudinger \cite{DT}. Moreover,
\begin{equation}\label{DT-emb}
W_0^{1,\Phi}(\Omega)  \stackrel{\hookrightarrow}{\mbox{\tiny cont}} L_{\Phi_*}(\Omega).
\end{equation}
\begin{remark} If $N \geq 3$ and $\phi(t) = 2$  then, by computing, we obtain for $t > 0$: $\Phi_*^{-1}(t) = t^{\frac{N-2}{2N}}$ and  $\Phi_{*}(t) = t^{\frac{2N}{N-2}}$. The function $\Phi_{*}$  plays the role of the critical Sobolev exponent in the case of Sobolev spaces.
\end{remark}
\begin{remark} The operator $\mbox{div} ( \phi(|\nabla u|) \nabla u)  = \Delta_{\Phi} u$ is referred to as the $\Phi$-Laplacian.
\end{remark}

\section{\bf Main Results}

\noindent Our main results are

\begin{theorem}\label{Teo_1}
        Assume $(\phi_1),(\phi_2)$ and $(\ref{cond_phi})$. Suppose there is a number  $ a\geq 0$ and a nonnegative function $b\in L^1(\Omega)$ such that
        \begin{equation}\label{cond_f}
            \displaystyle
           \begin{array}{l}
                |f(x,s)|\leq a\Phi_*(s)+b(x),~~ s\in \mbox{\bf R}~~ \mbox{a.e.}~~  x \in \Omega.
              \end{array}
           \end{equation}
\noindent  Assume also that  $F$ satisfies $(\ref{cond_F})$ and
        \begin{equation}\label{cond_teo}
            \inf_{v\in W_0^{1,\Phi}(\Omega),~\|v\|_\Phi=1}\left\{\int_\Omega \Phi(|\nabla v|)dx-\int_{\{v\neq 0\}}A_\infty(x)|v(x)|^\ell dx\right\}>0.
        \end{equation}
\noindent Then for  $h\in L_{{\Phi}}(\Omega)^{\prime}$, there is  $u\in W_0^{1,\Phi}(\Omega)$ satisfying  $(\ref{prob})$ in the sense of distributions.
       \end{theorem}

\begin{theorem}\label{Teor_2}
         Assume $(\phi_1),(\phi_2)$, $(\ref{cond_phi})$ and
        \begin{equation}\label{cond_f_fraco}
            \displaystyle
           |f(x,s)s|\leq a\Phi_*(s)+b(x)|s|,~~  s\in \mbox{\bf R}~~ \mbox{a.e.}~ x\in\Omega,
        \end{equation}
\noindent for some number $a\geq 0$, and a nonnegative $b \in L_{{{\Phi}_*}}(\Omega)^{\prime}$. Assume that  $F$ satisfies $(\ref{cond_F})$,  and  $(\ref{cond_teo})$ holds. Then for   $h\in L_{{\Phi}}(\Omega)^{\prime}$, there is a weak solution  $u\in  W_0^{1,\Phi}(\Omega)$ of  $(\ref{prob})$.
\end{theorem}
\begin{remark}  When $\phi(t) \equiv 1$  one has
\begin{itemize}
\item[(i)]~~  $W_0^{1,\Phi}(\Omega) = H_0^1(\Omega)$,
\item[(ii)]~~  equation $(\ref{prob})$ and condition $(\ref{cond_f_fraco})$ become respectively
\begin{equation}\label{Laplacian}
\displaystyle  -\Delta u = f(x,u) + h(x)~~ \mbox{in}~~
\Omega,
\end{equation}\label{critical}
\noindent and
$$
\displaystyle |f(x,s)|\leq a |s|^{2^* -1}+b(x),~~  \mbox{a.e.}~ x\in\Omega,~~  s\in   \mbox{\bf R},
$$
\item[(iii)]~~  condition $(\ref{cond_teo})$  becomes
\begin{equation}\label{cond_teo_DELTA}
            \inf_{v\in H_0^1(\Omega),~|v|_2=1}\left\{\int_\Omega |\nabla v|^2 dx-\int_{\{v\neq 0\}}A_\infty(x)v^2 dx\right\}>0,
        \end{equation}
\item[(iv)]~~ finding a weak solution of $(\ref{prob})$ means finding a weak solution  of $(\ref{Laplacian})$,

\item[(v)] when $A_\infty(x) \leq \alpha(x)~ \mbox{for some}~ \alpha \in L^{\infty}(\Omega)~ \mbox{with}~\alpha \leq \lambda_1~\mbox{in}~ \Omega,~ \mbox{and}~ \alpha < \lambda_1~\\
\mbox{on a subset of}~ \Omega~\mbox{with positive measure, where}~ \lambda_1~ \mbox{is the principal }\\
\mbox{eigenvalue of}~ (-\Delta, H_0^1(\Omega)),~ \mbox{then}~ (\ref{cond_teo_DELTA})$ holds, (\mbox{cf.}~ {\rm \cite{gonc})}.
\end{itemize}
\end{remark}
\noindent A classical result on integral equations which goes back to Hammerstein \cite{Hamm} shows that
\begin{equation}\label{Maw}
\displaystyle  -\Delta u = f(x,u) + h(x)~~ \mbox{in}~~
\Omega,~~ u = 0~~  \mbox{on}~~ \partial \Omega,
\end{equation}

\noindent is solvable provided $f(x,s)$  grows at most linearly in  s and a condition such as
$$
\displaystyle A_\infty(x) \leq \mu,~~  \mbox{a.e.}~~x \in \Omega~~  \mbox{for some}~~  \mu \in \mbox{\bf  R}.
$$
\noindent holds, with $\mu < \lambda_1$.  

\noindent  Mawhin, Willem \& Ward in \cite{Mawhin}  allowed subcritical
growth on $f(x, s)$ and a solution of (\ref{Maw}) was shown to exist under the additional condition
$$
  \begin{array}{l}
    \displaystyle A_\infty(x) \leq \alpha(x)~~ \mbox{a.e.}~x \in \Omega~~\mbox{for some}~ \alpha \in L^{\infty}(\Omega),~\mbox{with}\\
     \alpha \leq \lambda_1~\mbox{in}~ \Omega,~~\alpha < \lambda_1~\mbox{on a subset of}~\Omega~\mbox{with positive measure},   \\
\end{array}
 $$
\noindent   Goncalves in \cite{gonc}  allowed critical  growth condition on $f(x, s)$, obtaining solutions in the distribution sense under condition $(\ref{cond_teo_DELTA})$
which was introduced by Br\'ezis \& Oswald \cite{BO}.

\noindent In this paper we go back to the setting above regarding problem (\ref{prob}), this time  in the framework of Orlicz-Sobolev spaces.
 We refer the reader to the papers \cite{ Garcia,  hui-2, Fuk_2, Fuk_1, Fuk_0, Le, Gz2,  montenegro-2} and their references for nonlinear boundary value problems  on Orlicz-Sobolev spaces.

\noindent Problems envolving the $\Phi$-Laplacian operator appear in nonlinear elasticity, plasticity and generalized Newtonian fluids, see e. g.  \cite{Fuk_0}, \cite{Fuk_2} and their references.
\vskip.1cm

\noindent Consider the problem, (where the operator is an example of the general $\Delta_{\Phi}$ above),
\begin{equation}\label{app}
  -  \mbox{\rm div} \Big(\gamma\frac{(\sqrt{1+ |\nabla u|^2}-1)^{\gamma-1}}{\sqrt{1+ |\nabla u|^2}} \nabla u   \Big)   =  f(x,u)+ h~~ \mbox{in}~\Omega,~ u = 0~  \mbox{on}~ \partial \Omega,
\end{equation}
\noindent where $1 \leq \gamma<\infty$.

\begin{remark}
\noindent We shall use the notation  $\gamma^* = N \gamma/(N-\gamma)$ for  $\gamma  \in (1,N)$.
\end{remark}
\noindent The results below will be proved by applying theorems \ref{Teo_1} and \ref{Teor_2}.

\begin{theorem}\label{teo_app_1}
 Let $1 < \gamma < N$. Assume that
\begin{equation}\label{cond_fgamma}
            \displaystyle
           \begin{array}{l}
                |f(x,s)|\leq a|s|^{\gamma^{\star}}+b(x),~ s\in \mbox{\bf R}~ \mbox{a.e.}~  x \in \Omega,
              \end{array}
           \end{equation}
\noindent where  $ a\geq 0$ is some constant,   $b\in L^1(\Omega)$ is nonnegative  and
\begin{equation}\label{cond_Fgamma}
        \displaystyle
       F(x,s) \leq A|s|^\gamma + B(x),~  s \in \mbox{\bf R},~ \mbox{a.e.}~ x \in {\Omega},
         \end{equation}
\noindent for some constant $A \geq 0$ and some nonnegative function  $B \in L^1(\Omega)$. If  in addition,
\begin{equation}\label{cond_app_f}
            \inf_{v\in W_0^{1,\gamma},~|v|_{L^\gamma}=1}\left\{\int_\Omega \left(\sqrt{1+|\nabla v|^2}-1\right)^\gamma dx-\int_{\{v\neq 0\}}A_\infty(x)|v(x)|^\gamma dx\right\}>0,
   \end{equation}
\noindent then for each  $h\in L^{\frac{\gamma}{\gamma-1}}(\Omega)$  problem $(\ref{app})$ admits a solution  $u\in W^{1,\gamma}_0(\Omega)$, in the distribution sense.
\end{theorem}
\noindent The result  below is a variant of theorem \ref{teo_app_1} for the case of weak solutions.
\begin{theorem}\label{teo_app_2}
 Let $1 < \gamma < N$. Assume that
\begin{equation}\label{cond_fgamma_fraca}
            \displaystyle
           \begin{array}{l}
                |f(x,s)|\leq a|s|^{\gamma^{\star}-1}+b(x),~ s\in \mbox{\bf R}~ \mbox{a.e.}~  x \in \Omega,
              \end{array}
           \end{equation}
\noindent where  $ a\geq 0$ is some constant,   $b\in L^{\gamma^*}(\Omega)$ is nonnegative.
If  in addition, $(\ref{cond_Fgamma})$ and $(\ref{cond_app_f})$ hold then for each  $h\in L^{\frac{\gamma}{\gamma-1}}(\Omega)$  problem $(\ref{app})$ admits a weak solution  $u\in W^{1,\gamma}_0(\Omega)$.
\end{theorem}
\begin{remark} The function
$$
\phi(t)=pt^{p-2}\ln(1+t)+\frac{t^{p-1}}{t+1},~ t>0.
$$
\noindent where $1 < p < N - 1$, satisfies $(\phi_1),(\phi_2)$ and $(\ref{cond_phi})$. However, in this case $L_\Phi(\Omega)$ is not a Lebesgue space $L^q(\Omega)$. This follows by applying a result in {\rm \cite[pg 156]{Rao1}}.
\end{remark}
\begin{remark}
In our arguments, $C$ will  denote a positive $(cumulative)$ constant.
\end{remark}

\section{\bf Proofs of Theorems \ref{Teo_1} and \ref{Teor_2}}

\noindent  At first, we recall that $L_\Phi(\Omega)'=L_{\widetilde{\Phi}}(\Omega)$ and moreover,
$$
\langle h, u\rangle=\int_\Omega hudx,~ u\in L_{\Phi}(\Omega),
$$
\noindent  (cf. \cite[thm 8.19]{A}). Consider the energy functional asociated to  (\ref{prob}),
    $$
    I(u)=\int_\Omega\Phi(|\nabla u|)dx-\int_\Omega F(x,u)dx-\int_\Omega hudx,~u\in W_0^{1,\Phi}(\Omega).
    $$
\noindent It follows by using  $\Phi\in\triangle_2$ and  $(\ref{emb-ell})-(\ref{cond_F})$ that $I:W_0^{1,\Phi}(\Omega) \longrightarrow \mbox{\bf R}$ is defined.

\noindent Next we state and prove some technical lemmas.

\begin{lemma}\label{lema_sc}
Assume $(\phi_1)$, $(\phi_2)$, $(\ref{cond_phi})$  and $(\ref{cond_F})$. Then $I$ is weakly  lower semicontinuous, {\rm wlsc} for short.
\end{lemma}
\noindent  \proof Let  $(u_n)\subseteq W_0^{1,\Phi}(\Omega)$ such that $u_n\rightharpoonup u~~  \mbox{in}~ W_0^{1,\Phi}(\Omega)$. Then $u_n\rightarrow u$ in $L_\Phi(\Omega)$ and, by eventually passing to subsequences,
  $u_n \rightarrow u$ a.e. in $\Omega$~  and there is  $\theta_2\in L^{\ell}(\Omega)$ such that $|u_n|\leq \theta_2$ a.e. in $\Omega$. By $(\ref{cond_F})$  we have
$$
    F(x,u_n) \leq  \displaystyle A |\theta_2|^{\ell}+B(x).
      $$
\noindent   Since $F$ is a  Carath\'eodory function,
    $$F(x,u_n(x))\longrightarrow F(x,u(x))~ \mbox{a.e.}~ x\in\Omega.$$
\noindent  By Fatou's lemma,
    $$\limsup\int_\Omega F(x,u_n)dx\leq\int_\Omega F(x,u)dx.$$
\noindent Hence
   $$
    I(u)     \leq  \displaystyle\liminf\left\{\int_\Omega\Phi(|\nabla u_n|)dx-\int_\Omega F(x,u_n)dx-\int_\Omega hu_ndx\right\} = \displaystyle\liminf I(u_n)
$$
\noindent   showing that $I$ is {\rm wlsc}. $\hfill{\rule{2mm}{2mm}}$

 \begin{lemma}\label{Lema_sol}
        Assume $(\phi_1)$, $(\phi_2)$, $(\ref{cond_phi})$, $(\ref{cond_F})$ and $(\ref{cond_f})$.  
Let  $u\in  W_0^{1,\Phi}(\Omega)$ such that
  \begin{equation}\label{MIN-1}
    I(u)=\displaystyle \min_{v\in W_0^{1,\Phi}(\Omega)}I(v).
   \end{equation}
    Then $u$ satisfies $(\ref{prob})$ in the sense of distributions.
    \end{lemma}
\noindent \proof    Let $v\in{C}_0^\infty(\Omega)$ and $0<t<1$. Then $u+tv\in   W_0^{1,\Phi}(\Omega)$ and
   \begin{equation}
    \begin{array}{lll}
      0 & \leq & \displaystyle\frac{I(u+tv)-I(u)}{t}\label{minimum-pt} \\
        \\
        &   =  & \displaystyle\int_\Omega \Big[\frac{\Phi(|\nabla u+t\nabla v|)-\Phi(|\nabla u|)}{t}-
                 \frac{F(x,u+tv)-F(x,u)}{t}- hv \Big]dx
    \end{array}
    \end{equation}
\noindent   We claim that
   \begin{equation}\label{Phi-Gateaux}
    \displaystyle \lim_{t\rightarrow 0_+}\int_\Omega\frac{\Phi(|\nabla u+t\nabla v|)-\Phi(|\nabla u|)}{t}dx=\int_\Omega\phi(|\nabla u|)\nabla u
    \nabla vdx.
   \end{equation}
\noindent     Indeed, take a function  $\theta_t$ such that
    \begin{equation}\label{el1}
        \Phi(|\nabla u+t\nabla v|)-\Phi(|\nabla u|)=  \phi(\theta_t)\theta_t\left[|\nabla u+t\nabla v|-|\nabla u|\right]~\mbox{a.e. in}~\Omega
    \end{equation}
\noindent   and
         \begin{equation}\label{el2}
        \min\{|\nabla u+t\nabla v|,|\nabla u|\}  \leq   \theta_t \leq  \max\{|\nabla u+t\nabla v|,|\nabla u|\}~~ \mbox{a.e. in }~ \Omega.
\end{equation}
\noindent   By $(\ref{el2})$,  $\theta_t \rightarrow |\nabla u|$ a.e. in  $\Omega$ as $t\rightarrow 0_+$. We infer that
    \begin{equation}\label{el3}
        \displaystyle\lim_{t\rightarrow 0_+}\frac{\Phi(|\nabla u+t\nabla v|)-\Phi(|\nabla u|)}{t}=\phi(|\nabla u|)\nabla u\nabla v~~\mbox{a.e. in
        }~\Omega.
    \end{equation}
\noindent    By $(\ref{el1})$,  (\ref{el2}) and $(\phi_2)$ we get
   \begin{equation}\label{el4}
        \left|\displaystyle\frac{\Phi(|\nabla u+t\nabla v|)-\Phi(|\nabla u|)}{t}\right|\leq\phi(|\nabla u|+|\nabla v|)(|\nabla u|+|\nabla v|)|\nabla
        v|.
    \end{equation}
\noindent By  \cite[lemma A.2]{Fuk_1} one has  $\widetilde{\Phi}(t\phi(t))\leq \Phi(2t)$ for $t\in  \mbox{\bf R}$, so that
    $$\phi(|\nabla u|+|\nabla v|)(|\nabla u|+|\nabla v|)\in L_{\widetilde{\Phi}}(\Omega)$$
\noindent and by  the H\"older Inequality (cf. \cite{A}),
    $$
        \phi(|\nabla u|+|\nabla v|)(|\nabla u|+|\nabla v|) |\nabla v|\in L^1(\Omega).
   $$
\noindent  By $(\ref{el3})$, $(\ref{el4})$ and  Lebesgue's theorem, $(\ref{Phi-Gateaux})$ follows.
  \vskip.2cm

\noindent  We claim that
  \begin{equation}\label{LIM}
    \displaystyle \lim_{t\rightarrow 0_+}\int_\Omega\frac{F(x,u+tv)-F(x,u)}{t}dx=\int_\Omega f(x,u)vdx.
 \end{equation}

\noindent Indeed, take a function $\rho_{t}$ such that
    $$
    \frac{F(x,u+tv)-F(x,u)}{t}=f(x,\rho_t(x))v~\mbox{a.e. in}~\Omega
    $$
\noindent and
    \begin{equation}\label{ineq_above}
         \min\{u+tv,u\}  \leq \rho_t \leq \max\{u+tv, u\}~ \mbox{a.e. in }~ \Omega.
    \end{equation}
\noindent Using (\ref{DT-emb}) we infer that $\Phi_*(|u|+|v|)\in L^1(\Omega)$. Using (\ref{cond_f})  and (\ref{ineq_above}) we have
   $$
        \begin{array}{lll}
          |f(x,\rho_t)v|  & \leq & a\Phi_*(|u|+|v|)|v|+b|v| \\ \\
                        & \leq & \left(a\Phi_*(|u|+|v|)+b\right)|v|_\infty.
        \end{array}
    $$
\noindent  By  Lebesgue theorem,  $(\ref{LIM})$ follows. Passing to the limit in $(\ref{minimum-pt})$ we infer that $u$ is a distribution solution of (\ref{prob}). This proves lemma \ref{Lema_sol}. $\hfill{\rule{2mm}{2mm}}$
\vspace{1mm}

\noindent \proof ({\bf  of Theorem \ref{Teo_1}})~ At first we show that $I$ is coercive. Indeed, assume by the way of contradiction, that there is
$(u_n)\subseteq  W_0^{1,\Phi}(\Omega)$  such that
 $$
 \displaystyle\|\nabla u_n\|_\Phi\rightarrow \infty~~\mbox{and}~~~ I(u_n)\leq C.
 $$
\noindent  Using $(\ref{cond_F})$  and the H\"older Inequality  we have
    \begin{equation}\label{egt1}
        \displaystyle \int_\Omega \Phi(|\nabla u_n|)dx \leq \displaystyle A\int_\Omega|u_n|^\ell dx+2\|h\|_{\widetilde{\Phi}}\|u_n\|_\Phi+C
      \end{equation}
\noindent  We claim that $\displaystyle \int_\Omega |u_n|^\ell dx\rightarrow\infty$. Indeed, assume on the contrary, that
$$
    \displaystyle\int_\Omega |u_n|^\ell dx\leq C.
    $$
\noindent By   (\ref{egt1}) and  Poincar\'e's Inequality (cf. \cite{Gz1}),
$$
\displaystyle \int_\Omega \Phi(|\nabla u_n|)dx \leq  C(1 +\|\nabla u_n\|_\Phi),
$$
\noindent which is impossible because  by  \cite[lemma 3.14]{Gz1},
$$
\displaystyle\frac{\displaystyle\int_\Omega\Phi(|\nabla u_n|)dx}{\|\nabla u_n\|_\Phi}\rightarrow\infty,
$$
\noindent showing that $\displaystyle \int_\Omega |u_n|^\ell dx\rightarrow\infty$. We infer, using $(\ref{emb-ell})$  that
    $\|u_n\|_{\Phi}\rightarrow\infty$.

\noindent By (\ref{egt1}) and lemma \ref{lema_naru} in the Appendix, we have
\begin{equation}\label{egt2}
          \|\nabla u_n\|_\Phi^\ell \leq \displaystyle \int_\Omega \Phi(|\nabla u_n|)dx \leq C\|u_n\|_{\Phi}^\ell+2\|h\|_{\widetilde{\Phi}}\|u_n\|_\Phi+C.
    \end{equation}
\noindent Dividing in $(\ref{egt2})$  by $\|u_n\|_\Phi^\ell$ we get
    $$\displaystyle\|\nabla v_n\|_\Phi^\ell\leq C+\frac{2\|h\|_{\widetilde{\Phi}}}{\|u_n\|_\Phi^{\ell-1}}+\frac{C}{\|u_n\|_\Phi^\ell},$$
 \noindent where $\displaystyle v_n=\frac{u_n}{\|u_n\|_\Phi}$. It follows that  $(\|\nabla v_n\|_\Phi)$ is bounded.
\noindent Passing to a subsequence, we have,
\begin{itemize}
      \item $v_n\rightharpoonup v$ in $W_0^{1,\Phi}(\Omega)$~~ and~~ $\displaystyle\int_\Omega\Phi(|\nabla v|)dx\leq\liminf\int_\Omega\Phi(|\nabla v_n|)dx$,
      \item $v_n\rightarrow v$ in $L_\Phi(\Omega)$~ and there is~ $\theta_3\in L^{\ell}(\Omega)~ \mbox{such that}~|v_n|\leq\theta_3$~ a.e. in ~$\Omega.$
  \end{itemize}
\noindent  We claim that
    \begin{equation}\label{Af1}
      \displaystyle\limsup\int_\Omega\frac{F(x,\|u_n\|_\Phi v_n)}{\|u_n\|_\Phi^\ell}dx\leq\int_{\{v\neq 0\}}A_\infty(x)|v(x)|^\ell dx.
    \end{equation}
\noindent  Indeed, it follows using  (\ref{cond_F}) that
    \begin{equation}\label{tcd1}
     \displaystyle\frac{F(x,\|u_n\|_\Phi v_n(x))}{\|u_n\|_\Phi^\ell}  \leq  \displaystyle A\theta_3^\ell(x)+B(x).
       \end{equation}
\noindent   In addition,
    $$\displaystyle \limsup\frac{F(x,\|u_n\|_\Phi v_n)}{\|u_n\|_\Phi^\ell}\leq\limsup\frac{F(x,\|u_n\|_\Phi
    v_n)}{\|u_n\|_\Phi^\ell|v_n(x)|^\ell}|v_n(x)|^\ell\chi_{\{v_n\neq 0\}}$$
\noindent  and hence
    \begin{equation}\label{tcd2}
          \displaystyle\limsup\frac{F(x,\|u_n\|_\Phi v_n)}{(\|u_n\|_\Phi v_n(x))^\ell}|v_n(x)|^\ell\chi_{\{v_n\neq 0\}}  \leq
          A_\infty(x) |v(x)|^\ell,~  v(x)\neq 0.
    \end{equation}
\noindent   By $(\ref{tcd1})$, $(\ref{tcd2})$ and  Fatou's Lemma, $(\ref{Af1})$ follows.  Using again the fact that $\Phi$ is  convex and continuous, lemma \ref{lema_naru}, $I(u_n) \leq C$ and  $(\ref{Af1})$ it follows   that
     $$
    \begin{array}{lll}
       \displaystyle\int_\Omega\Phi(|\nabla v|)dx & \leq & \displaystyle \liminf\int_\Omega\Phi(|\nabla v_n|)dx\\ \\
       & \leq & \displaystyle \liminf\frac{1}{\|\ u_n\|_\Phi^\ell}\int_\Omega\Phi(|\nabla u_n|)dx\\ \\
       & \leq & \displaystyle \liminf\left\{\int_\Omega\left[\frac{F(x,\|u_n\|_\Phi
       v_n)}{\|u_n\|_\Phi^\ell}+\frac{2\|h\|_{\widetilde{\Phi}}}{\|u_n\|_\Phi^{\ell-1}}\right]dx+\frac{C}{\|u_n\|_\Phi^\ell}\right\}\\ \\
       & \leq & \displaystyle \limsup\int_\Omega\frac{F(x,\|u_n\|_\Phi v_n)}{\|u_n\|_\Phi^\ell}dx\\ \\
       & \leq & \displaystyle\int_{\{v\neq 0\}}  A_\infty(x) |v(x)|^\ell dx,
    \end{array}
    $$
\noindent   which contradicts (\ref{cond_teo}). Therefore,  $I$  is coercive. By lemma \ref{lema_sc},  there is $u\in W_0^{1,\Phi}(\Omega)$ satisfying $(\ref{MIN-1})$ and by lemma \ref{Lema_sol},  $u$ satisfies (\ref{prob}) in the sense of distributions. This proves theorem \ref{Teo_1}.  $\hfill{\rule{2mm}{2mm}}$
\vspace{1mm}

\noindent  The lemma below is needed in order to prove theorem \ref{Teor_2}.
    \begin{lemma}\label{Lema_sol_fraco}
       Assume $(\phi_1)$, $(\phi_2)$,  $(\ref{cond_phi}),  (\ref{cond_F})$ and $(\ref{cond_f_fraco})$.
    If $u\in  W_0^{1,\Phi}(\Omega)$ satisfies $(\ref{MIN-1})$  then $u$ is a weak solution of $(\ref{prob})$.
    \end{lemma}
 \noindent    The proof is similar to that of lemma \ref{Lema_sol}.  In the present case one must show
$(\ref{Phi-Gateaux})$ and $(\ref{LIM})$ for $v \in W_0^{1,\Phi}(\Omega)$ which follows basically the same lines as in the case $v \in C_0^{\infty}$.

\noindent  As earlier there is  a function $\rho_{t}(x)$ such that\\
$$
    \frac{F(x,u+tv)-F(x,u)}{t}=f(x,\rho_t(x))v,~\mbox{a.e. in}~\Omega
    $$
\noindent  and
    $$
         \min\{u+tv,u\}  \leq \rho_t \leq \max\{u+tv, u\}~ \mbox{a.e. in }~\Omega.
   $$
\noindent  Using (\ref{cond_f_fraco}) and the fact that $t\mapsto\frac{\Phi_*(t)}{t}$ is increasing (cf. \cite{A}) we get to
    $$
        |f(x,\rho_t)v|   \leq  a\Phi_*(|u|+|v|)+b|v|.
    $$
\noindent  Since $\Phi_*(|u|+|v|)\in L^1(\Omega)$ we get $f(x,\rho_t)v\in L^1(\Omega)$.
Applying Lebesgue's theorem we get $(\ref{LIM})$. As in the proof of lemma \ref{Lema_sol} we infer that $u$ is a weak solution of ($\ref{prob}$). This proves  lemma \ref{Lema_sol_fraco}.    $\hfill{\rule{2mm}{2mm}}$
\vspace{1mm}
   
\noindent \proof ({\bf  of Theorem \ref{Teor_2}})  As in the proof of theorem \ref{Teo_1} one shows that  $I$ is coercive. Since by lemma \ref{lema_sc}, $I$ is {\rm wlsci}, there is  $u\in  W_0^{1,\Phi}(\Omega)$ such that satisfying $(\ref{MIN-1})$. By lemma \ref{Lema_sol_fraco}, $u$ is a weak solution of  $(\ref{prob})$. This proves theorem \ref{Teor_2}. $\hfill{\rule{2mm}{2mm}}$

\section {\bf Proofs of Theorems \ref{teo_app_1} and \ref{teo_app_2}}

\noindent  We shall need some preliminary results. Set
    $$
    \phi(t)=\gamma\frac{(\sqrt{1+t^2}-1)^{\gamma-1}}{\sqrt{1+t^2}}, ~t\geq 0.
    $$
\begin{lemma}\label{equivalenceSob}
       The function  $\phi$ satisfies $(\phi_1),(\phi_2)$,~ and
    \begin{equation}\label{newDelta2}
    \gamma\leq\frac{t^2\phi(t)}{\Phi(t)}\leq 2\gamma,~~t>0.
    \end{equation}
 \noindent  In addition,   when $1 < \gamma < N$,
 $$
 L_\Phi(\Omega)=L^\gamma(\Omega),\\
 $$
\begin{equation}\label{embedding}
 |u|_\Phi\leq |u|_\gamma,~u\in L_\Phi(\Omega),
\end{equation}
\noindent  and
 $$  W_0^{1,\Phi}(\Omega)=W^{1,\gamma}_0(\Omega).$$
\end{lemma}

\noindent \proof Of course $\phi\in{C}(0,\infty)$. By a direct computation we infer that for $t > 0$,
    $$ \lim_{t\rightarrow 0} t\phi(t)=0,~~   \displaystyle \lim_{t\rightarrow \infty} t\phi(t)=\infty,~~  (t\phi(t))' > 0~ \mbox{and}~ \gamma \leq \frac{t^2\phi(t)}{\Phi(t)}   \leq   2\gamma, $$
 \noindent showing $(\phi_1)$, $(\phi_2)$ and  $(\ref{newDelta2})$. To prove that $L_\Phi(\Omega)=L^\gamma(\Omega)$, we point out that
$$
    \Phi(t)\leq t^\gamma,~~t\geq 0~~ \mbox{and}~~  \Phi(t) \geq \frac{1}{2^\gamma}t^\gamma,~~t \geq 2.
    $$
\noindent By a result in  \cite[p 156]{Rao1},  $L_\Phi(\Omega)=L^\gamma(\Omega)$. As a consequence, $W_0^{1,\Phi}(\Omega) =W^{1,\gamma}_0(\Omega)$,
\vskip.1cm

\noindent In order to show $(\ref{embedding})$, take $u\in L^\gamma(\Omega)$ and $k>0$ and notice that
$$
\displaystyle \int_\Omega \frac{|u(x)|^\gamma}{k^\gamma}\leq 1~~\mbox{iff} ~~|u|_\gamma\leq k.
$$
\noindent Since  $\Phi(t)\leq t^\gamma~\mbox{for}~t\geq 0$ we have
$$
\displaystyle \int_\Omega\Phi\left(\frac{u}{k}\right)dx\leq\frac{|u|_\gamma^\gamma}{k^\gamma}.
$$
\noindent Setting  $k=|u|_\gamma$ we get $(\ref{embedding})$. This  proves  lemma \ref{equivalenceSob}.

\begin{proposition}\label{pro_app_2}
   Assume $(\ref{cond_app_f})$. Then $(\ref{cond_teo})$ holds.
  \end{proposition}
\noindent \proof  By $(\ref{cond_app_f})$ and $W_0^{1,\Phi}(\Omega)=W^{1,\gamma}_0(\Omega)$,
$$
\displaystyle \inf_{u \in W_0^{1,\Phi}(\Omega),~|u|_{L^\gamma}=1}\Big\{\int_\Omega \big(\sqrt{1+|\nabla u|^2}-1\big)^\gamma dx-\int_{\{u\neq 0\}}A_\infty(x)|u(x)|^\gamma dx\Big\}>0,
$$
\noindent Recalling that $\Phi(t)=(\sqrt{1+t^2}-1)^{\gamma}$, there  is $\delta > 0$ such that
$$
\delta \leq\displaystyle \int_\Omega\Phi\left(\frac{|\nabla u|}{|u|_\gamma}\right)dx-\frac{1}{|u|_\gamma}\int_\Omega A_\infty(x)|u|^\gamma dx,~~u\in W_0^{1,\Phi}(\Omega).
$$
\noindent Using the convexity of $\Phi$ and $|u|_\Phi\leq |u|_\gamma$  we have
$$
\delta   \leq\displaystyle \int_\Omega\Phi\left(|\nabla u|\right)dx-\int_\Omega A_\infty(x)|u|^\gamma dx,~~u\in W_0^{1,\Phi}(\Omega),~ |u|_\Phi  = 1.
$$
\noindent This proves proposition \ref{pro_app_2}.

\begin{proposition}\label{pro_app_3}  Assume $1 < \gamma < N$ and $(\ref{cond_fgamma})$. Then  $ (\ref{cond_f})$ holds.
   \end{proposition}

\noindent \proof  Set  $\Phi(t)=t^{\gamma}$. By lemma \ref{lema_naru_*} (in the Appendix), we have for $t \geq 1$,
    $$ \Phi_*(1)~t^{\gamma^*}\leq \Phi_*(t).$$
\noindent  Using the inequality above and $(\ref{cond_fgamma})$ we get $ (\ref{cond_f})$. This proves  proposition \ref{pro_app_3}.
\vspace{2mm}

\noindent   \proof({\bf of Theorem \ref{teo_app_1}}) As a consequence of  propositions \ref{pro_app_2} and \ref{pro_app_3},
theorem \ref{Teo_1} applies ending the proof of theorem \ref{teo_app_1}. $\hfill{\rule{2mm}{2mm}}$
\vspace{1mm}

\noindent  \proof ({\bf  of Theorem \ref{teo_app_2} }) Similarly to the proof of the theorem above it suffices to apply theorem \ref{Teor_2}.$\hfill{\rule{2mm}{2mm}}$

 \section {\bf Appendix}

\noindent  We refer the reader to  \cite{Fuk_1} for the  lemmas below whose proofs are elementary.

\begin{lemma}\label{lema_naru}
       Assume that  $\phi$ satisfies  $(\phi_1)-(\phi_2)$ and $(\ref{cond_phi})$.
        Set
         $$
         \zeta_0(t)=\min\{t^\ell,t^m\},~~~ \zeta_1(t)=\max\{t^\ell,t^m\},~~ t\geq 0.
        $$
\noindent  Then  $\Phi$ satisfies
       $$
            \zeta_0(t)\Phi(\rho)\leq\Phi(\rho t)\leq \zeta_1(t)\Phi(\rho),~~ \rho, t> 0,
        $$
$$
\zeta_0(\|u\|_{\Phi})\leq\int_\Omega\Phi(u)dx\leq \zeta_1(\|u\|_{\Phi}),~ u\in L_{\Phi}(\Omega).
 $$
\end{lemma}
\begin{lemma}\label{lema_naru_*}
    Assume that  $\phi$ satisfies $(\phi_1)-(\phi_2)$ and $(\ref{cond_phi})$.  Set
    $$
    \zeta_2(t)=\min\{t^{\ell^*},t^{m^*}\},~~ \zeta_3(t)=\max\{t^{\ell^*},t^{m^*}\},~~  t\geq 0.
    $$
\noindent    Then
        $$
            \ell^*\leq\frac{t^2\Phi'_*(t)}{\Phi_*(t)}\leq m^*,~t>0,
       $$
        $$
            \zeta_2(t)\Phi_*(\rho)\leq\Phi_*(\rho t)\leq \zeta_3(t)\Phi_*(\rho),~~ \rho, t> 0,
       $$
       $$
            \zeta_2(\|u\|_{\Phi_{*}})\leq\int_\Omega\Phi_{*}(u)dx\leq \zeta_3(\|u\|_{\Phi_*}),~ u\in L_{\Phi_*}(\Omega).
        $$
\end{lemma}

\enddocument